\documentclass[12pt]{amsart}
\usepackage{amsfonts,amsmath,amsthm,amscd,amssymb,latexsym,cite,verbatim,texdraw,floatflt,
caption2}
\RequirePackage{hyperref}

\usepackage[a4paper]{geometry}

\theoremstyle
{plain}

\begin{document}

\title{Coarse  spaces, ultrafilters and dynamical  systems }

\author{ Igor Protasov}

\maketitle
\vskip 5pt

{\bf Abstract.}  
For a coarse space   $(X, \mathcal{E})$,  $X^\sharp$ denotes the set of all unbounded ultrafilters on $X$ 
endowed with the parallelity relation: $p||q$  if there exists 
 $E \in \mathcal{E} $   such that  $   E[P]\in q  $ for each $P\in p$.
If    $(X, \mathcal{E})$ is finitary then there exists a group  $G $   of permutations of $X$   such that the coarse structure  $\mathcal{E}$  has the base $\{\{ (x,gx): x\in X$,  $g\in F\}:  F\in  [G]^{<\omega},   \  id \in F \}.$
We survey and analyze interplays between  $(X, \mathcal{E})$,  $X^\sharp$  and the dynamical system $(G, X^\sharp)$.

\vspace{6 mm}


\vspace{3 mm}


\vspace{10 mm}

The dynamical $\check{S}$varc-Milnor Theorem and Gromov Theorem arose at the dawn of {\it Geometric Group Theory.}
In both cases, a group or a pair of groups act on some locally compact spaces, see
 [22, Chapter 1].
The Gromov coupling criterion was transformed into the powerful tool in coarse equivalences (see references in 
\cite{b23}), 
 however some natural questions on the coarse  equivalence of groups  need more delicate combinatorial  technique, see 
\cite{b4}. 

In this paper, we describe and survey the dynamical approach to coarse spaces originated in the algebra of the  Stone-$\check{C}$ech compactification. 
We identify the Stone-$\check{C}$ech compactification $\beta G$ of  a discrete group $G$ with the set of all ultrafilters on $G$. 
The left regular action $G$ on $G$ gives  rise to the action of $G$  on  $\beta G$  by
$(g,p)\mapsto  gp $, $gp=\{ gP: P\in p \}$. 
In turn on, the dynamical system $(G,\beta G)$ induces on $\beta G$ the structure of a right topological semigroup. 
The product $pq$ of ultrafilters $p,q$ is defined by $A\in pq$ if and only if $\{ g\in G: g^{-1} A\in q\}\in p.$ 
The semigroup $\beta G$ has very rich algebraic structure and the plenty of combinatorial applications, see nice paper
\cite{b5},
capital book 
\cite{b6} 
 or booklet 
\cite{b9}. 

Let $(X, \mathcal{E})$ be a coarse space. 
We denote by $X^\sharp$ the set of all ultrafilters $p$ on $X$ such that each member $P\in p$ is unbounded in $(X, \mathcal{E})$. 
Then we define the parallelity equivalence $||$ on $X^\sharp$ by $p||q$ if and only if there exists $E\in \mathcal{E}$ such that $E[P]\in q$ for each $P\in p$. 
For $p\in X^\sharp$, the orbit
$    \overline{ \overline{p}} = \{ q\in X^\sharp : q || p \} $
 looks like a smile apart of some hidden cat. 
This cat appears if
$(X, \mathcal{E})$ 
  is finitary. 
By Theorem 3.1, there exists a group $G$ of permutations of $X$ such that
$ \mathcal{E}$ 
  has the base $\{\{ (x, gx): x\in F \}: F\in [G]^{< \omega} $,  $id \in F \}.$
In this case,
$X^\sharp = X^\ast$,  $X^\ast =\beta X \setminus X$  and $ \overline{ \overline{p}} =Gp$. 
But even 
$(X, \mathcal{E})$ 
 is not finitary, $X^\sharp$ contains some counterpart of the kernel of a dynamical system, see Theorem 2.3.

\vspace{6 mm}

Our goal is to clarify interplays between 
$(X, \mathcal{E})$,
  $X^\sharp$
   and the dynamical system $(G, X^\ast)$ in order to understand the dynamical nature of some extremal coarse spaces,  in particular, tight, discrete and indiscrete.

\section{ Coarse   spaces}

Given a set $X$, a family $\mathcal{E}$ of subsets of $X\times X$ is called a {\it coarse structure }  on $X$  if
\vskip 7pt

\begin{itemize}
\item{}   each $E\in \mathcal{E}$  contains the diagonal  $\bigtriangleup _{X}$,
$\bigtriangleup _{X}= \{(x,x)\in X: x\in X\}$;
\vskip 5pt

\item{}  if  $E$, $E^{\prime} \in \mathcal{E}$ then $E\circ E^{\prime}\in\mathcal{E}$ and
$E^{-1}\in \mathcal{E}$,   where    $E\circ E^{\prime}=\{(x,y): \exists z((x,z) \in  E,  \   \ (z, y)\in E^{\prime})\}$,   $E^{-1}=\{(y,x): (x,y)\in E\}$;
\vskip 5pt

\item{} if $E\in\mathcal{E}$ and $\bigtriangleup_{X}\subseteq E^{\prime}\subseteq E  $   then
$E^{\prime}\in \mathcal{E}$;
\vskip 5pt

\item{}   $\bigcup  \mathcal{E}= X\times X $.

\end{itemize}
\vskip 7pt

A subfamily $\mathcal{E}^{\prime} \subseteq \mathcal{E}$  is called a
{\it base} for $\mathcal{E}$  if,
 for every $E\in \mathcal{E}$, there exists
  $E^{\prime}\in \mathcal{E}^{\prime}$  such  that
  $E\subseteq E ^{\prime}$.
For $x\in X$,  $A\subseteq  X$  and
$E\in \mathcal{E}$, we denote
$$E[x]= \{y\in X: (x,y) \in E\},
 \   E [A] = \bigcup_{a\in A}   \    
  E[a], \   \     
   E_A   [x]  = E[x]\cap A$$
 and say that  $E[x]$
  and $E[A]$
   are {\it balls of radius $E$
   around} $x$  and $A$.

The pair $(X,\mathcal{E})$ is called a {\it coarse space}  
 \cite{b22} 
 or a ballean
\cite{b16 }, 
\cite{b21}. 

For a coarse   space $(X,\mathcal{E})$, a  subset $B \subseteq X$   is called {\it bounded} if $B \subseteq  E[x]$  for some
$E\in \mathcal{E}$ and $x\in X$.
The family $\mathcal{B}_{(X, \mathcal{E})}$
 of all bounded subsets of  $(X, \mathcal{E})$
is called the {\it bornology} of $(X, \mathcal{E})$.

A coarse  space  $(X,\mathcal{E})$ is called {\it  finitary, } if for each  $E\in   \mathcal{E} $ there exists a natural number $n$ such that $|E[x]|< n$  for each $x\in X$.

We classify subsets of a coarse space  $(X,\mathcal{E})$ by their size.
 A  subset $A$ of  $X$ is called
\vskip 7pt
\begin{itemize}

\item{}   {\it  large} if  $E[A]=X$   for some
$E\in \mathcal{E};$
\vskip 10pt

\item{}   {\it  small} if  $L \backslash A$ is large for
each large subset $L$;
\vskip 10pt

\item{}   {\it  thick} if,  for each  $E\in \mathcal{E}$,
there exists 
 $a\in A$  such that 
$E[a]   \subseteq A;$
\vskip 10pt

\item{}   {\it  prethick} if  $E[A]$  is thick for some 
   $E\in \mathcal{E}$;

\item{} {\it thin (or discrete)} if, for each
 $E\in \mathcal{E}$, 
 there
exists a bounded subset 
$B$ of $X$ such that
 $E_A[a]= \{a\}$ for each $a\in A\backslash B$.

\end{itemize}

\vskip 7pt

For finitary coarse spaces, the dynamical unification of above definitions  will  be given  in Section 3. 
\vskip 10pt

Following \cite{b17},  we  say that two subets $A, B$  of $X$ are 
\vskip 7pt
\begin{itemize}

\item{}   {\it close} (write $A \delta  B$)  if there  exists $E\in \mathcal{E}$  such that,
 $A \subseteq  E[B]$,  $B \subseteq  E[B]$;

\item{}     {\it  linked }  (write $A\lambda   B$) if either $A, B$ are  bounded or there exist unbounded  subsets $A^{\prime} \subseteq  A$,  $B^{\prime} \subseteq  B$ such that  $A^{\prime}\delta   B^{\prime}$.
\end{itemize}
\vskip 7pt

We  say that a coarse space $E\in \mathcal{E}$   is 

\vskip 7pt

\begin{itemize}

\item{}  {\it $\delta$-tight}  if any two unbounded subsets of $X$ are close;

 \item{}  {\it $\lambda$- tight}   if any two unbounded subsets of $X$ are linked;

 \item{}  {\it indiscrete}   if   $E\in \mathcal{E}$  has no 
 unbounded discrete subsets;

 \item{}  {\it ultradiscrete}   if   $\{ X \backslash   B: B\in  \mathcal{B}_{(X, \mathcal{E})}\} $
  is an
 unltrafilter.

\end{itemize}

We note that  $\lambda$-tight spaces appeared in 
\cite{b2} 
  under the name 
{\it   utranormal}, $\delta$-tight subsets are called 
{\it extremely normal} in 
\cite{b14}
 and {\it hypernormal} in  
 \cite{b1}.

An unbounded  coarse space  is called {\it maximal} if it is bounded
in every stronger coarse structure. 
By [18, Theorem 3.1], 
 every maximal coarse space is $\delta$-tight.
A ballean  $(X,\mathcal{E})$  is 
$\delta$-tight if and only if  every subset of $X$  is large. 
If a $\lambda$-tight space is not indiscrete  then it contains an   
 ultradiscrete subspace 
[14, Theorem 2.2], 
  so every finitary  $\lambda$-tight  space is indiscrete.


\section{Ultrafilters}

Let  $X$ be a discrete space and let $\beta X$ denotes the  $Stone-\check{C}ech \  \ compactification$ of  $X$.
 We take the points of $\beta X$ to be the ultrafilter on $X$, with the points of $X$ identified with the principal ultrafilters, 
so $X^\ast = \beta X\setminus X$ is the set of all free ultrafilters.  The topology of  $\beta X$  can be defined by stating that 
the sets of the form
 $\bar{A} =\{ p\in \beta X: A\in p \}$,
 where $A$ is a subset of $X$,  are base for the open  sets. 
The universal property of $\beta X$ states that every mapping  
$f: X\longrightarrow Y$, where $Y$ is a compact Hausdorff space, can be extended to the continuous mapping
$f^\beta : \beta X\longrightarrow X$.

Given a coarse space $(X, \mathcal{E})$, we endow $X$ with the discrete topology and denote by $X^\sharp$ the set of all  ultrafilters $p$ on $X$ such that each member $P\in p$ is unbounded. Clearly, $X^\sharp$ is the closed subset of $X^\ast$ and $X^\sharp= X^\ast$ if
 $(X,\mathcal{E})$
  is finitary.

Following 
\cite{b10}, 
 we say that two ultrafilters $p, q\in X^\sharp$ are {\it parallel} (and write $p || q$) if there exists  $E\in \mathcal{E}$ such that $E[P]\in q$ for each $P\in p$. 
By [10, Lemma 4.1, 1], 
$||$ is an equivalence on $X^\sharp$. We denote by $\sim$ the minimal (by inclusion) closed (in  $X^\sharp \times  X^\sharp$) equivalence on $X^\sharp$ such that $|| \subseteq \sim$. 
The quotient
$\nu (X,\mathcal{E})$ of 
 $X^\sharp$ by $ \sim$
  is called the Higson corona of   $(X,\mathcal{E})$.
For $p\in X^\sharp$, we denote 
$$    \overline{ \overline{p}} = \{ q\in X^\sharp : q || p \} ,  \   \   \Breve{p}  = \{ q\in p :  q \sim p \} .$$

A function
$f: (X,\mathcal{E}) \longrightarrow   \mathbb{R}$
  is called {\it slowly oscillating} if, for every  
$E\in \mathcal{E}$ 
and $\epsilon >0$,  there exists a bounded subset $B$ of $X$ such that $diam  f(E[x])<\epsilon$ for each $x\in X\setminus B$.

We recall \cite{b10} 
 that a coarse space  $ (X,\mathcal{E})$ is {\it normal} if any two asymptotically disjoint subsets $A, B$ of $X$
have disjoint asymptotic neighbourhoods. Two subsets $A,  B$ of $X$  are called {\it asymptotically disjoint}  if $E[A] \cap E[B]$ is bounded for each $E\in \mathcal{E}$.
A subset $U$ of $X$ is called an {\it asymptotic neighbourhood} of a subset $A$ if $E[A]\setminus U$ is bounded for each 
 $E\in \mathcal{E}$.
By [10,  Theorem 2.2], 
 $ (X,\mathcal{E})$   is normal if and only if, 
for any two disjoint and asymptotically  disjoint subsets $A, B$  of $X$, there exists a slowly oscillating function $f: X\longrightarrow  [0,1]$ such that $f|_ A =0$,  $f|_B =1$.

By   [11, Proposition 1], 
  $p\sim q$  if and only if 
$h^\beta (p) = h^\beta (q)$
for every slowly oscillating function 
$h: (X,\mathcal{E})\longrightarrow    [0,1]. $

By [4,  Theorem 7], 
 $ (X,\mathcal{E})$   is normal if and only if
$\sim \  = \  c l  \  ||$.

By [17, Theorem 9 and Corollary 10], 
 if 
$\lambda _{ (X,\mathcal{E})} =  \lambda _{ (X,\mathcal{E} ^\prime )} $  then
 Higson coronas  of   $(X,\mathcal{E})$  and  $(X,\mathcal{E} ^\prime )$  coincide and if    $(X,\mathcal{E})$  is normal then  $(X,\mathcal{E} ^\prime )$   is normal.

By [21, Theorem 2.1.1] 
 a coarse space   $ (X,\mathcal{E})$  is metrizable if $ \mathcal{E}$  has a countable base.
If   $\lambda _{(X,  \mathcal{E})}=\lambda _{(X,  \mathcal{E^{\prime}})}$
and  $ (X,\mathcal{E})$ is metrizable then $(X,  \mathcal{E^{\prime}})$ needs not to be metrizable 
 [17,  Theorem 3]. 

\vskip 10pt

{\bf Question 2.1 [17].}  {\it  Let
$\delta _{(X,  \mathcal{E})}=$
  $\delta _{(X,  \mathcal{E^{\prime}})}$ and 
$(X, \mathcal{E})$ is  metrizable. Is 
 $(X, \mathcal{E^{\prime}})$  metrizable?}

\vskip 10pt

If the answer to Question 2.1 would be positive then  
 $ \mathcal{E} = \mathcal{E}^\prime$, 

 Let 
$(X, \mathcal{E})$ be a coarse space. We say that a subset  $S$  of  $X^\sharp$ is {\it invariant} if 
$\overline{\overline{p}} \subseteq S$ for each $p\in S$. 
Every non-empty closed invariant subset of $X^\sharp$ contains a minimal by inclusion closed invariant subset. 
We denote 

$ \  \   \  K(X^\sharp) = \bigcup \{ M: M  \  \  $ is minimal closed invariant subset of $  \  X^\sharp \}.$

\vskip 10pt

{\bf Theorem 2.2.} {\it For $ p\in  X^\sharp $,  $cl \overline{\overline{p}}$  is a minimal closed invariant subset if and only if, for every  $P\in p$, there exists 
 $E\in \mathcal{E}$ such that $ \overline{\overline{p}} \in (E[P])^\sharp$. 
\vskip 7pt

Proof.  } Apply arguments proving this statement for metric spaces
 [13,  Theorem 3.1].
$ \ \  \  \Box$ 
\vskip 10pt

{\bf Theorem 2.3.} {\it For $ q\in  X^\sharp $,  $q \in cl K   (X^\sharp)$ if and only if each subset  $Q\in q$ is prethick.

\vskip 7pt

Proof.  } Apply arguments proving   Theorem 3.2  in 
 \cite{b13}. 
$ \ \  \  \Box$ 

{\bf Theorem 2.4.} {\it Let  $p,q$ be  ultrafilters  from $X^\sharp$ such that 
$\overline{\overline{p}}$,  $\overline{\overline{q}}$ are countable and  
$cl \overline{\overline{p}}\cap  cl\overline{\overline{q}} \neq \emptyset.$ 
Then either $cl \overline{\overline{p}}\subseteq  cl\overline{\overline{q}} $  or 
$cl \overline{\overline{q}}\subseteq  cl\overline{\overline{p}}. $

\vskip 7pt

Proof.  } Apply arguments proving   Theorem 3.4  in 
 \cite{b13}. 
$ \ \  \  \Box$ 
\vskip 10pt


\section{Dynamical systems}

By a {\it dynamical system} we mean a pair  $(G, T)$, where $T$ is a compact space, $G$ is a group of homeomorphisms of $G$.

The following two theorems make a bridge between 
coarse spaces and dynamical systems.
 For usage of  Theorem 3.1 in corona constructions see 
 \cite{b3}. 

Let $G$ be a 
transitive
 group  of permutations of a set  $X$. 
We denote by $X_G$  the  set  $X$  endowed with the coarse structure with  the base. 
$$   \{\{ (x, gx): g\in F \} : F \in  [G]^{<\omega},  \   id \in  F\}. $$

\vskip 10pt

{\bf Theorem 3.1.} {\it 
For every finitary coarse space   $(X, \mathcal{E})$, there exists  a group $G$  of permutations of $X$ such  that  
 $(X, \mathcal{E})= X_G .$

\vskip 7pt

Proof.  }
Theorem 1 in
 \cite{b12},
  for more   general  results  see 
\cite{b8}, 
\cite{b15}.
$ \ \  \  \Box$ 

\vskip 10pt


{\bf Theorem 3.2.} {\it     If  
$(X, \mathcal{E})$,  $(X, \mathcal{E}^{\prime}  )$
are finitary coarse   spaces and  
$||_{(X, \mathcal{E})} =$  $ ||_{(X, \mathcal{E}^{\prime}  )}$
then  $ \mathcal{E} = \mathcal{E}^{\prime}$.

\vskip 7pt

Proof.  }
Theorem 15 in 
\cite{b17}.
$ \ \  \  \Box$ 

\vskip 10pt

If   $(X, \mathcal{E})=X_G$, 
we say that  $X_G$ is the $G$-{\it realization }  of 
$(X, \mathcal{E})$. Each
$G$-realization  of $(X, \mathcal{E})$
   defines the dynamical system 
$(G, X^\ast)$
  with the action 
 $(g,p) \mapsto  gp$,
$gp= \{ gP: P\in p \}$.
We note that  
 $ \overline{\overline{p}} = Gp$ for each  $p\in X^\ast$.
Since the parallelity $||$ is defined by means of entourages,  
  the partition of   $ X^\ast$ into $G$-orbits  does not depend on $G$-realizations  of $(X, \mathcal{E})$.
It follows that if some property formulated in terms of $G$-orbits of $ (G, X^\ast)$  is proved for some $G$-realization of 
 $(X, \mathcal{E})$  then it holds for any  $G$-realization. 
Moreover, by Theorem 3.2, every finitary coarse structure can be uniquely recognized by the set of orbits in $X^\ast$.

\vskip 10pt

Given a finitary coarse space  $(X, \mathcal{E})$, its  $G$-realization of 
 $ X_ G$, a subset  $A\subseteq  X$ and $p\in X^\ast$, we define  the {\it $p$-companion} of $A$  by
$$\triangle _p (A)= A^\ast \cap Gp. $$

\vskip 10pt

{\bf Theorem 3.3.} { \it For a subset $A$ of  $(X, \mathcal{E})$,  the following statements hold

\vskip 10pt

$(1)$   $A$  is large iff  $\triangle _p (A)\neq \emptyset  $ for each $p\in X^\ast$; 
\vskip 7pt

$(2)$	$A$  is thick iff  $\triangle _p (A)= G p  $ for  some $p\in X^\ast$; 
\vskip 7pt

$(3)$	$A$  is thin  iff  $|\triangle _p (A)| \leq 1$ for each $p\in X^\ast$;    
\vskip 7pt

\vskip 7pt

Proof.  }  Theorem 3.1 and 3.2   in
 \cite{b20} .
$ \ \  \  \Box$ 

\vskip 10pt

We recall that a  dynamical  system $(G, T)$ is 

\vskip 7pt
\begin{itemize}

\item{}  {\it minimal}  if each orbit $Gx$ is dense in $T$;\vskip 7pt

\item{}  {\it topologically transitive} if some orbit  $Gx$  is dense in $T$.

\end{itemize}

\vskip 10pt

For a dynamical system $(G, T)$,  $ker (G, T)$ denotes the closure of the union of all minimal  closed $G$-invariant  subsets of $T$. 
  Theorem 2.3  describes  explicitely the kernel of the dynamical  system $(G,  X^\ast)$  of  $X_G$.

\vskip 10pt

{\bf Theorem 3.4.} { \it Let  $(X, \mathcal{E})$ be a finitary coarse space 
 and   $(X, \mathcal{E})= X_G$.  Then $(X, \mathcal{E})$ is $\delta$-tight if and only if the dynamical  system $(G, X^\ast)$ is minimal.

\vskip 7pt

Proof.  }  Apply  Theorem 3.3(1). $ \ \  \  \Box$ 

\vskip 10pt

{\bf Theorem 3.5.} { \it Let  $(X, \mathcal{E})$ be a finitary coarse space 
 and   $(X, \mathcal{E})= X_G$.
  Then  the following statements are equivalent 

\vskip 10pt

$(1)$   $(X, \mathcal{E})$ is $\lambda$-tight ; 
\vskip 7pt

$(2)$	 for any infinite subset   $A, B$ of  $X$,  there exist  
 $p\in X^\ast$  and  $g\in G$  such that $A\in p$,  $B\in gp$; 
\vskip 7pt

$(3)$	for any  family   $\{A_n : n\in \omega \}$  of infinite subsets of  $X$,   
 there exists  
 $p\in X^\ast$  such that 
$A_n ^\ast   \cap Gp\neq \emptyset  $ for each $n\in \omega$.
\vskip 10pt

Proof.  } $(1)\Longrightarrow (2)$.  Since $A, B$ are linked, there exist $A^\prime \subseteq A $, $B^\prime \subseteq  B$ and $H\in [G]^{<\omega}$  such that  $A^\prime \subseteq H B^\prime $. We take 
$p\in X ^\ast$  such  $A^\prime \in p$. 
Then $B^\prime  \in  h^{-1}p $ for some $h\in H$.
\vskip 7pt 

$(2)\Longrightarrow (3)$.  We choose inductively a sequence  $(g_n)_{n\in\omega}$ in $G$ and
a sequence $(C_n) _{n\in\omega}  $
  of subsets of $G$  such  that  $C_n\subseteq A_n$,  $g_n C_n\subseteq A_{n+1}$, $C_{n+1}\subseteq g_n C_n$.
Let  $ h_n = g_n g_{n-1} \dots  g_0$.  Then
$$ A_0 \cap  h_0 ^{-1} A_1 \cap  \dots \cap  h_n ^{-1} A_{n+1}  \neq\emptyset $$
for each $n\in \omega. $ 
We take an arbitrary ultrafilter $p\in X^\sharp$ such that 
$ A_0 \cap  h_0 ^{-1}  \cap  \dots \cap  h_n ^{-1} A_{n+1}  \in p $ for each $n\in\omega$. 
Then  $Gp\cap  A_n ^\ast  \neq\emptyset $  for each $n\in\omega$.

 \vskip 7pt 

$(3)\Longrightarrow (1)$.  Evident.
  $ \ \  \  \Box$

\vskip 10pt

{\bf Corollary 3.6.}  {\it If  $(X, \mathcal{E})= X_G$     and the dynamical  system  $(G, X^\ast)$ is topologically  transitive then $(X, \mathcal{E})$ is $\lambda$-tight.}

\vskip 10pt

{\bf Remark 3.7.}  Does there exist a group $G$ of permutations of $\omega$ such that   $\omega_G$ is $\lambda$-tight and $(G, \omega^\ast)$ is not topologically transitive? 
This question can not be answered in ZFC without additional assumptions. Yes, if $\frak{t}< \frak{c}$ and No  if 
 $\frak{t}= \frak{c}$, 
 see [1, Theorems 5.2 and 5.3]. 

\vskip 10pt

{\bf Theorem 3.8.} { \it Let $K$ be a closed nowhere dense subset of $\omega^\ast$. 
Then there exists a transitive group $G$ of permutations of $\omega$ such that $ker (G, \omega^\ast) = K$ and the orbit  $Gp$ is dense in $\omega^\ast$ for each $p\notin K$.

\vskip 10pt

Proof.  } We take a filter $\phi$ on $\omega$  such  that  $K=\overline{\phi}$ and $\overline{\phi}$   has the base $\{\overline{A}:  A\in\phi\}$. 
We denote by $G$ the group of all permutations $g$ of $\omega$ such that there exists $A_g$ such that $g(x)=x$ for each $x\in A_g$.
Clearly, $G$ is transitive. 

If $q\in K$ then $g(q)=q$ for each $g\in G$ so $K\subseteq ker (G, \omega^\ast)$.

We fix   $p\in \omega^\ast\setminus K$ and take an arbitrary $q\in\omega^\ast$, $p\neq q$.  Let $P\in p$, $Q\in q$ and $P\cap Q =\emptyset$. Since $K$ is nowhere dense, there exists $A\in \phi$ such that $P\setminus A\in p$ and $Q\setminus A$ is infinite. 
By the definition of $G$, there exists $g\in G$ such that $g(P\setminus A)= Q\setminus A$. Hence, $Gp$ is dense in $\omega^\ast$ and $K= ker (G, \omega^\ast)$.
$ \ \  \  \Box$ 

\vskip 10pt

{\bf Corollary 3.9.}
{\it There are $2^{\frak{c}}$ $\lambda$-tight finitary coarse spaces on $\omega$ which are not $\delta$-tight.

\vskip 10pt

Proof.  } In light of Corollary 3.6, it suffices to  notice that there are $2^{\frak{c}}$  free ultrafilters on $\omega$.
$ \ \  \  \Box$ 

\vskip 10pt

Each orbit of a dynamical system from the proof of Theorem 3.8 is either dense or a singleton. We construct a topologically transitive $(G, \omega^\ast)$ having an infinite discrete orbit.

\vskip 10pt

{\bf Example 3.10. }
We partition $\omega$  into infinite subsets $\{ W_n : n\in  \mathbb{Z}\}$,   fix a bijection  $f_n : W_n  \longrightarrow  W_{n+1}$ and denote by $f$  a  bijection of  $\omega$  such that    $f|_{W_n} = f_n$. 

For each $n\in \mathbb{Z}$,  we pick  $p_n \in \omega^\ast$  such that  $W_n\in  p_n$  and denote by  $S$   the set of all  permutations  $g$  such that,  for each 
$g(x)=x$, $x \in W_n$.

We take  the group  $G$  of  permutations  generated by  $S\cup \{f \}$.
Then  $G p_0 = \{ p_n : n\in \mathbb{Z}\}$ 
and $Gp_0$ is discrete.

If  $p\in W_0 ^\ast$   and  $p\neq  p_0$ then  $Gp$ is  dense in $\omega^\ast$ so  $(G, \omega^\ast)$ is topologically transitive.

\vskip 10pt


{\bf Remark 3.11. }
If
 $(X, \mathcal{E})$
 is a finitary coarse space, 
 $(X, \mathcal{E})=X_G$
 then, by Theorem 2.3, every infinite subset of 
 $(X, \mathcal{E})$
 is prethick if and only if $ker (G, X^\ast) = X^\ast$.

If every infinite subset of a finitary coarse space
 $(X, \mathcal{E})$
  is thick then
 $(X, \mathcal{E})$
 is discrete.
 Indeed, if
 $(X, \mathcal{E})$
 is not discrete then, by Theorem 3.3(3), there exists 
$q\in X^\ast$ and $g\in G$ such that $gq\neq q$. 
We take $Q\in q$ such that $gQ\cap Q=\emptyset$. 
It follows that $Q$ is not thick.

\vskip 10pt

We say that a subset $A$ of $X_G$  is 

\vskip 7pt
\begin{itemize}

\item{}  {\it sparse}  if    $\triangle_p (A) $  is finite   for each     $p\in G^\ast$;\vskip 7pt

\item{}  {\it scattered} if,   for each infinite subset        $Y$ of  $A$   there exists $p\in Y^\ast$ such that  $\triangle_p (Y) $ is finite.

\end{itemize}

\vskip 10pt

{\bf Theorem 3.12.} {\it Let $(X, \mathcal{E})$ be a finitary coarse space and $(X, \mathcal{E})=X_G$. Then
the following statements are equivalent: 
\vskip 10pt

$(1)$   $(X, \mathcal{E})$  is indiscrete; 
\vskip 7pt

$(2)$	for every infinite subset   $A$  of   $ X$,  there exist  $p\in X^\ast$  and $g\in G$ such that $A\in p$,  $A\in g p$
and $p\neq gp$;
    
\vskip 7pt

$(3)$	every infinite subset $A$ of $X$   is not sparse.

\vskip 10pt

Proof. } The equivalence of $(1)$ and $(2)$  follows from  Theorem 3.3(3),  $(3)\Longrightarrow (1)$
is evident. 

To show $(2)\Longrightarrow (3)$, we choose a sequence  $(g_n)_{n\in\omega}$ in $G$ and sequence 
$(A_n)_{n\in\omega}$ of subsets of $A$ such that 
$$A_{n+1} \subset A_{n}, \   \   g_n  A_{n+1}\cap  A_{n+1} = \emptyset, \  \  n\in\omega   $$
and choose $p\in X^\ast$ such that $A_n \in p$  for each $n\in\omega$.
Then $Gp\cap A^\ast$ is infinite and  $A$ is not sparse.  
 $ \ \  \  \Box$ 

\vskip 10pt

{\bf Theorem 3.13.} {\it  A subset $A$  of    $X_G$  is scattered  if  and  only if $Gp$ is  discrete for each  $p\in A^\ast$.

\vskip 10pt

Proof. }  Theorem 5.4  in   
\cite{b20} . 
$ \ \  \  \Box$ 
\vskip 10pt

Let $G$ be a group of permutations of a set $X$. 
Let $(g_n)_{n\in \omega}$ 
 be a sequence in $G$ and let
$(x_n)_{n\in \omega}$
 be a sequence in $X$ such that
\vskip 7pt

$(1)$   $  \    \{ g_0 ^{\epsilon_0} \dots  g_n ^{\epsilon_n} \ x_n : \epsilon _i \in \{0,1\} \}
\cap
\{ g_0 ^{\epsilon_0} \dots  g_m ^{\epsilon_m} \  x_m : \epsilon _i \in \{0,1\} \}=\emptyset$
 for all distinct $n,m\in \omega$;

\vskip 7pt

$(2)$   $  \    \{ g_0 ^{\epsilon_0} \dots  g_n ^{\epsilon_n} \ x_n : \epsilon _i \in \{0,1\} \}
|=2 ^ {n+1}$
 for every  $n\in \omega$.
\vskip 10pt

Following [20], we say that a subset $Y$ of $X$ is a {\it piece-wise shifted $FP$-set} if there exist  
$ (g_n)_{ n\in \omega}$,  $ (x_n)_{ n\in \omega}$
satisfying (1), (2) and such that 

$$ Y =   \{ g_0 ^{\epsilon_0} \dots  g_n ^{\epsilon_n} \ x_n : \epsilon _i \in \{0,1\} \}, \ 
n\in \omega \}.$$

\vskip 10pt

{\bf Theorem 3.14.} {\it A subset $A$ of $X_G$ is scattered if and only if $A$ does not contain piece-wise shifted $FP$-sets.
\vskip 10pt

Proof.} Theorem 4.4 in [20].
$  \   \Box$
\vskip 10pt

{\bf Theorem 3.15.} {\it
Let   $(X, \mathcal{E})$   be a finitary indiscrete space, 
$(X, \mathcal{E})=X_G$.
Then there exists $p\in X^\ast$ such that the orbit $G_p$ is not discrete.
\vskip 10pt

Proof.} 
We may suppose that $G$ consists of all permutations $g$ of $X$ such that $(x, gx)\in E$ for some $E\in  \mathcal{E}$ and all $x\in X$.

In light of Theorem 3.13 and Theorem 3.14, it suffices to find a piece-wise shifted $FP$-set in $X$ defined by some sequence
$(g_n )_ {n\in\omega}$
  in $G$ and some sequence
$(x_n )_ {n\in\omega}$
  in $X$.

Since $X$ is not discrete, there are an infinite subset $A_0 \subset X $ and an involution $g_0\in G$ such that  $A_0 \cap g_0 A_0 = \emptyset$ and $g_0 x =x$ for each $x\in X\setminus (A_0 \cup g A_0)$. 
Pick $x_0\in A_0$.

Suppose that $A_0, \dots , A_n$, $g_0, \dots , g_n$ and
 $x_0, \dots , x_n$
 have been chosen. 
Since $A_n$ is not discrete, we can find  an infinite subset $A_{n+1}$ of $A_n$ and an involution  $g_{n+1} \in G$  such that $g_{n+1} A_{n+1} \subset A_n$,  
$A_{n+1}\cap  g_{n+1} A_{n+1} = \emptyset$
and $g_{n+1} x = x $ for each $x\in X\setminus (A_{n+1} \cup g_{n+1} A_{n+1})$.
Pick $x_{n+1}\in  A_{n+1}$.

After $\omega$ steps, we get the desired sequences  $(g_n)_{n\in \omega}$,
$(x_n)_{n\in \omega}$.  $ \ \Box$

\vskip 10pt

{\bf Example 3.16. } We show that
an indiscrete
  space needs not to be $\lambda$-tight. 
To this end, we take the set $X$  of all rational  number on  $[0,1]$, denote  by  $G$ the group of all homeomorphisms of  $X$ and consider the finitary space $X_G$.
Let  $A=\{ a_n: n\in \omega \}$,  $B=\{ b_n: n\in \omega \}$ be subsets  of $X$  such that $(a_n)_{n\in\omega}$
converges to $0$ and  $(b_n)_{n\in\omega}$   converges to some irrational number. 
Then $A, B$ are not linked, so  $X_G$ is not  $\lambda$-tight. 
On the other hand,
let $A$ be an infinite subset of $X$. 
We take two disjoint sequences $(a_n)_{n\in\omega}$,   $(b_n)_{n\in\omega}$  in $A$ which converge to some point $x\in X$. Then there is a homeomorphism $g$ of $X$ such that $g(a_n)=b_n$, $n\in \omega$. 
Applying Theorem 3.12, we see that $X_G$ is indiscrete.       
\vskip 10pt

As the results, we have got the following  line
\vskip 7pt

$\delta$-tight $\Longrightarrow$ topologically transitive  $\Longrightarrow$   $\lambda$-tight  $\Longrightarrow$  indiscrete,  

in which the first and the third arrow can not be reversed, and the second arrow can be reversed but only under some assumptions additional to ZFC.

\vskip 10pt

{\bf Question 3.17.}  {\it  Is a dynamical system $(G, \omega^\ast)$ minimal provided that $ker (G, \omega^\ast)=\omega^\ast$ and $(G, \omega^\ast)$ is topologically transitive?}

 \vskip 10pt      

The Higson corona of  every $\lambda$-tight space is a singleton. 

\vskip 10pt

 The following example suggested by Taras Banakh shows that the Higson corona of finitary indiscrete space
needs not to be a singleton.
\vskip 10pt

{\bf Example 3.18.}   Let $(X_1, \mathcal{E}_1)$, $(X_2, \mathcal{E}_2)$  be infinite indiscrete finitary spaces. We endow the union $X$ of $X_1$  and $X_2$
with the smallest coarse structure $\mathcal{E}$  such that the restrictions of  $\mathcal{E}$ to $X_1$  and $X_2$  coincide with $\mathcal{E}_1$  and $\mathcal{E}_2$.
Then the Higson corona of $(X, \mathcal{E})$  is not a singleton because the function $f$ defined by $f(x)=0$,  $x \in X_1$  and  $f(x)=1$,  $x \in X_2$  is slowly oscillating. 

\vskip 10pt

A subset $A$  of a coarse space  $(X, \mathcal{E})$   is  called  {\it $n$-thin (or $n$-discrete)}, 
$n\in \mathbb{N}$ if for each $E\in \mathcal{E}$ there  exists a bounded subset $B$  of 
$(X,\mathcal{E})$ such that $|E_A [a]|\leq n$  for each  $a\in A\setminus B$. 
Every $n$-thin metrizable  coarse space  can be partitioned into $\leq n$  thin subsets 
\cite{b7},
 but the Bergman's  construction from 
\cite{b19} 
gives a finitary $n$-thin space which can not be 
partitioned into $\leq n$  thin subsets.

\vskip 10pt

{\bf 
Theorem 3.19.} {\it  There exists a group $G$ of permutations of  $\omega$ such that the
coarse space $\omega_G $ is 2-discrete but $G$ can not be finitely partitioned into 
discrete subsets.
\vskip 10pt

Proof.} Theorem 6.1 in [1]. $ \  \Box$

\vskip 10pt

{\bf Acknowledgements. } I thank Taras Banakh for discussions around $G$-realizations of finitary spaces


\vskip 15pt

CONTACT INFORMATION
\vskip 15pt

I.~Protasov: \\
Faculty of Computer Science and Cybernetics  \\
        Kyiv University  \\
         Academic Glushkov pr. 4d  \\
         03680 Kyiv, Ukraine \\ i.v.protasov@gmail.com

\end{document}